\documentstyle[12pt]{article}
\newcommand{\const}{\mathop{\rm const}\limits}

\newcommand{\Var}{\mathop{\rm Var}\limits}

\newcommand{\argmax}{\mathop{\rm argmax}\limits}

\begin{document}

 \begin{center}

{\bf GENERALIZED BERNSTEIN-TYPE APPROXIMATION    \\

\vspace{3mm}

OF CONTINUOUS FUNCTIONS. } \\

\vspace{4mm}

{\sc Eugene Ostrovsky, Leonid Sirota}\\

\vspace{3mm}

 Bar-Ilan University,  59200, Ramat Gan, ISRAEL; \\
 e-mail: eugostrovsky@list.ru \\
 e-mail: sirota3@bezeqint.net \\

\vspace{4mm}

     {\sc Abstract.}

\end{center}

 \ We  derive in this short article
the non-asymptotical non-uniform sharp error estimation for the Bernstein's type approximation of
continuous function based on the modern probabilistic apparatus.\par

\vspace{4mm}

\par

\vspace{3mm}

\section{ Introduction. Notations. Statement of problem.}

\vspace{3mm}

 \ Let $ I = (a,b) $ or $ I = [a,b) $ or $ I = [a,b], a < b  $ be  a finite or infinite segment on the real axis, $  f: I \to R $
be {\it  uniformly } continuous  bounded function, $  \{\xi_i \} = \xi_i(x), \ i = 1,2,\ldots,n; \ \xi = \xi_1 = \xi(x),
 \ x \in I  $ be a {\it family } of  independent identically distributed (i., i.d.)  random variables (r.v.) which values in the set $  I, $
and  whose distribution dependent on the parameter $  x $ such that

$$
{\bf E} \xi_i = x, \ \Var \xi_i = \sigma^2(x). \eqno(1.1)
$$

 \ {\it It will be presumed that}  $  0 < \sigma(x) < \infty. $  Evidently, if $ - \infty < a < b < \infty,  $ then
$ \sigma(x) \le 0.5 (b - a).  $ \par

\ Denote

$$
S_n = n^{-1} \sum_{i=1}^n \xi_i, \hspace{5mm} S_n^o = n^{-1} \sum_{i=1}^n \xi_i - x. \eqno(1.2)
$$

 \ Let  us introduce the following {\it sequence of approximated linear operators}

$$
A_n[f](x) \stackrel{def}{=} {\bf E} f(S_n) = {\bf E} f( x + S_n^o), \eqno(1.3)
$$
which are in turn some generalization of the classical Bernstein's operators. See also \cite{Achiezer1}, \cite{Chen1}, \cite{DeVore1},
\cite{Floater1}, \cite{Gonska1}, \cite{Lorentz1}, \cite{Stancu1}, \cite{Timan1}. \par
 \ For instance,

$$
B_n[f](x) = \sum_{m=0}^n  {n \choose m } \ f \left(\frac{m}{n} \right) \ x^m \ (1 - x)^{n-m}
$$
be the ordinary Bernstein's polynomial of degree $  \ n, \ $  \cite{Bernstein1}, see also  \cite{Bojanic1}, \cite{Bojanic2},
\cite{Li2},  \cite{Lorentz1},  \cite{Math'e}, \cite{Palatanea1}, \cite{Palatanea2}, \cite{Telyakovskii1}. \par

\vspace{3mm}

 \ Here
$$
I = [0,1], \ a = 0, b = 1, \  {\bf P}(\xi = 1) = x, \
$$

$$
{\bf P}(\xi = 0) = 1 - x, \ \sigma(x) = [ x(1 - x)]^{1/2}.
$$

\vspace{3mm}

 \ Another ("Poisson") example. Here

$$
a = 0, \ b = \infty, \ I = [0, \infty), \ {\bf P }(\xi = k) = e^{-x} \ x^k /k!, \ k = 0,1,2,\ldots,
$$
so that $ {\bf E} \xi = x = \Var \xi, \ \sigma(x) = \sqrt{x},   $

$$
S_n[f](x) = e^{-nx} \sum_{k=0}^{\infty} \frac{(n x)^k}{k!} \ f \left( \frac{k}{n}  \right),
$$
see e.g. \cite{Chlodovsky1}, \cite{Chlodovsky2},  \cite{Mursalem1}, \cite{Szasz1}.

\vspace{4mm}

\ {\bf  We intend to investigate in this report the non-asymptotical
(quality and quantity) error of the uniform approximation properties of introduced operators:  } \par

\vspace{4mm}

$$
\Delta_n = \Delta_n[f] \stackrel{def}{=} \sup_{x \in I} | \ A_n[f](x) - f(x) \ |. \eqno(1.4)
$$

 \ Note in addition  that the standard deviation (weigh) $  \sigma(x) $ may has a form

$$
\sigma(x) = C \ x^{\alpha} \ (1 - x)^{\beta}, \alpha, \beta = \const \ge 0, \ I = [0,1];
$$
so that $ \sigma(x)  $  is Jacobi weight, see  \cite{Glazyrina1}. \par

\vspace{3mm}

\section{ Auxiliary apparatus and notations.}

\vspace{3mm}

 \hspace{4mm} We will use the following concrete  Ditzian-Totik  \cite{Ditzian1}   modulus of continuity:

$$
\omega_{\sigma}[f](\delta) \stackrel{def}{=} \sup_{|h| \le \delta} \sup_{x \in I} | \ f(x + h \sigma(x)) - f(x) \  |, \eqno(2.1)
$$
where we agree to take $  f(x)  = f(b), $ if $ x > b, $ and $  f(x) = f(a)  $ in the case $ x < a. $ \par
\ As ordinary modulus of continuity, the Ditzian-Totik modulus
is continuous monotonically increasing function equal to zero at the origin. \par

\ Recall that the {\it tail function} $ T_{\eta} (u) $ for each (not necessary to be
 non - negative)  r.v. $  \eta $ is defined as follows:

$$
T_{\eta} (u) = {\bf P}( |\eta| > u), \ u \ge 0. \eqno(2.2)
$$

\vspace{4mm}

 {\bf Definition 2.1.} Let $  \eta  $ be a {\it  centered } r.v. with finite variance, and let $  \eta(i) $ be independent
copies of the r.v. $ \ \eta \ $ defined perhaps  on some sufficiently rich probability space.  The following numerical function
$  Q(u) = Q_{\eta}(u), \ u \ge 0 $ will be named as {\it an absolute tail function,} write  $ ATF = ATF(\eta) $ for the r.v. $  \eta: $

$$
Q(u) = Q_{\eta}(u) \stackrel{def}{=} \sup_n {\bf P} \left(  \ n^{-1/2} \ \left| \sum_{i=1}^n  \eta(i) \ \right| > u  \right),
\ u \ge 0. \eqno(2.3)
$$

 \ This function was introduced and partially investigated by S.N.Bernstein;
we will bring for the estimation of these function  the most advanced probabilistic methods. \par
\vspace{4mm}

{\bf Definition 2.2.} The function $  u \to T(u), \ u \ge 0 $ is said to be {\it tail  function,} briefly, TF,
if it is right continuous, monotonically non - decreasing, $  T(0) = 1, $  and $  T(\infty) = 0. $ \par
 \ Obviously, every such a function is tail function for some random variable. \par

\vspace{4mm}

 \ It is clear that $ Q_{\eta}(u)  $ is tail function, as long as $  0 < \Var \eta < \infty. $ It follows from the classical CLT
that  for  some positive constant $  K = K(\eta)  $

$$
Q_{\eta}(u) \ge  \exp  \left(  - K \ u^2  \right), \ u \ge 1.
$$

 \ Note in addition that if the mean zero r.v. $  \eta $ is bilateral bounded: \\
 $ - \infty < a \le \eta \le b < \infty, $
then

$$
Q_{\eta}(u) \le  \exp  \left(  - K_2 \ u^2  \right), \ u \ge 1,
$$
Hoeffding's inequality.  Moreover, this estimate there holds iff the mean zero r.v. $  \eta $ is subgaussian,
see \cite{Buldygin1} -  \cite{Buldygin3},  \cite{Kahane1}.  \par

\vspace{3mm}

 \ Many examples of estimation of ATF function may be found in \cite{Kozatchenko1};  \cite{Ostrovsky1}, chapter 1,
sections 1.6.,  chapter 2. \par

\vspace{3mm}

 \ Further, let us consider  the centered and normed variable (variate)

$$
\zeta_i = \zeta_i(x) := \frac{\xi_i - x}{\sigma(x)}; \ \zeta  = \zeta(x) := \zeta_1. \eqno(2.4)
$$

\ Define the logarithm of  moment generating function for the r.v. $  \zeta:  $

$$
\phi(\lambda) = \phi_{\zeta}(\lambda) \stackrel{def}{=} \max_{\pm} \sup_{x \in I} \ln {\bf E} \exp ( \pm \lambda \ \zeta),
$$
if of course  the r.v. $ \zeta $ satisfies the famous Kramer's  condition, which is equal in turn the following implication

$$
\exists \lambda_0 \in (0,\infty], \ \Rightarrow  \forall \lambda:  \ |\lambda|  < \lambda_0  \ \Rightarrow
\phi(\lambda)  < \infty.
$$
 \ This function is also even, convex and generated the so-called Banach space $ B(\phi) $ consisting on the special centered
random variables, see \cite{Kozatchenko1}, \cite{Ostrovsky1}, chapters 1,2.\par

 \ Denote

$$
\nu(\lambda) = \sup_n [ \ n \ \phi(\lambda/\sqrt{n}) \ ],
$$

$$
\nu^*(u) = \sup_{\lambda} (\lambda u - \nu(\lambda)),
$$
the Young-Fenchel transform for the function $ \nu(\cdot). $  Both the introduced functions $ \nu(\cdot), \ \nu^*(\cdot) $
are correct definite, even and convex.  The following Fenchel-Moraux identity play a very important role in the theory of
random variables with exponential decreasing tails of distributions: $  \nu^{**}(u) = \nu(u).  $ \par

\ It is known, see  \cite{Ostrovsky1}, chapter 1, that

$$
Q_{\zeta}(u) \le 2 \ \exp ( - \nu^*(u)  ), \ u \ge 0. \eqno(2.5)
$$

\vspace{3mm}

\section{ Main result.}

\vspace{3mm}

\ {\bf Theorem 3.1.}

$$
\Delta_n[f] \le \int_0^{\infty}  \omega_{\sigma} [f] \left( \frac{z}{\sqrt{n}}  \right) \  \left| d Q_{\zeta}(z) \right|. \eqno(3.0)
$$

\vspace{3mm}

\ {\bf Proof.} We have using the direct definition of Ditzian-Totik modulus of continuity

$$
\Delta_n[f] \le \sup_x {\bf E} \left| f( x +  \zeta_n \sigma(x) /\sqrt{n}) - f(x)  \right| \le
$$

$$
 {\bf E} \sup_x \left| f( x +  \zeta_n \sigma(x) /\sqrt{n}) - f(x)  \right| \le
 {\bf E} \omega_{\sigma}[f] \left( |\zeta_n| \ /\sqrt{n} \right) =
$$

$$
 \int_0^{\infty} \omega_{\sigma}[f] \left( z /\sqrt{n} \right) \ | d T_{\zeta_n}(z) |. \eqno(3.1)
$$

\vspace{3mm}

 {\bf Lemma 3.1.} Let $ g = g(z), \ z \ge 0 $ be a
 continuous monotonically increasing function equal to zero at the origin. If $ \xi, \ \eta  $ are two non-negative
r.v. such that

$$
T_{\xi}(z) \le T_{\eta}(z), \ z \ge 0,
$$
then

$$
{\bf E} g(\xi) \le {\bf E} g(\eta). \eqno(3.2)
$$

\vspace{3mm}

{\bf Proof} of lemma 3.1. We can and will assume without loss of generality that all the functions $  g(x), T_{\xi}(z) $
and $ T_{\eta}(z)  $ are continuous differentiable. We deduce by means of integration by parts:

$$
{\bf E} g(\xi) = - \int_0^{\infty} g(x) \ d T_{\xi}(x) = - g(x) \ T_{\xi}(x) \  /_0^{\infty} +
\int_0^{\infty} T_{\xi}(x) \ g'(x) \ dx  =
$$

$$
\int_0^{\infty} T_{\xi}(x) \ g'(x) \ dx  \le \int_0^{\infty} T_{\eta}(x) \ g'(x) \ dx = -\int_0^{\infty} g(x) \ d T_{\eta}(x)  =
{\bf E} g(\eta).
$$

\vspace{3mm}

 \ There is an another proof. Namely, we can realize both the r.v. $ \xi, \ \eta $ on at the same probability space, say
$   [0,1], $ so that

$$
\eta = T_{\eta}^{-1}(\tau),  \hspace{6mm}   \xi = T_{\xi}^{-1}(\tau),
$$
where the r.v. $ \ \tau \ $ has an uniform distribution on the unit interval $  [0,1], $ if for definiteness both the tail
functions $ T_{\eta}(\cdot) $  and $  T_{\eta}(\cdot) $ are continuous and strictly decreasing. \par
 \ Therefore $  \xi \le \eta  $ almost  everywhere in this realization and following  $ {\bf E} g(\xi)  \le  {\bf E} g(\eta) $
under arbitrary realization.  \par

\vspace{3mm}

 \ It is no hard to finish the proof proof of theorem 3.1. Since   $ T_{\zeta_n}(z) \le  Q_{\zeta}(z),  $
we conclude on the basis of lemma 3.1

$$
\Delta_n[f] \le
 \int_0^{\infty} \omega_{\sigma}[f] \left( z /\sqrt{n} \right) \ | d Q_{\zeta}(z) |,
$$
Q.E.D. \par

\vspace{3mm}

{\bf Remark 3.1.} It follows immediately from the proposition of theorem 3.1 by virtue of Lebesgue dominated convergence theorem
that under formulated above conditions

$$
\lim_{n \to \infty} \Delta_n[f] = 0,
$$
as long as $ \ \omega_{\sigma}[f](\delta) \le 2 \sup_x |f(x)|.  $\par

\vspace{3mm}

 {\bf Example 3.1.} Suppose that the (centered normed) variable $  \zeta $ has a following tail function

$$
T_{\zeta}(u) \le \exp \left( - u^p  \right), \ u \ge 0 \eqno(3.2)
$$
for some constant  $  p > 0. $ Denote $  q = q(p) = \min(p,2).  $ It is known, see
\cite{Kozatchenko1}, \cite{Ostrovsky1}, chapters 1,2 that

$$
Q_{\zeta}(u)  \le \exp \left( - K(p) \ u^q  \right), \ u \ge 0, \ K(p)  = \const \in (0, \infty),  \eqno(3.3)
$$
and the last estimate is essentially non - improvable. \par
 \ We get relying on the theorem 3.1

$$
\Delta_n[f]  \le K(p) \int_0^{\infty}  \omega_{\sigma} [f] \left( \frac{z}{\sqrt{n}}  \right) \
 z^{q-1} \exp \left( - K(p) \ z^q   \right) \ dz. \eqno(3.4)
$$

\vspace{3mm}

\section{Some examples.}

\vspace{3mm}

 \hspace{4mm} We will consider in this section some examples in order to make sure the result of theorem 3.1. Note at first that
the case of the classical Bernstein's approximation in this spirit was considered in  \cite{Ostrovsky201}. \par

\vspace{3mm}

{\bf Definition 4.1.} The (continuous) function $  f: I \to R $ belongs by  definition to the  H\"older -
Ditzian - Totik class,  write $  f \in HDT(\sigma, \alpha), $ iff

$$
\omega_{\sigma}[f](\delta) \le H \cdot \delta^{\alpha}, \ \delta \ge 0, \eqno(4.1)
$$
for some constants $  0 \le H < \infty, \ \alpha \in (0,1]. $ \par

 \ We will understood as a capacity of the value $  H  $ in (4.1) its minimal value, namely

$$
H = H_{\alpha,\sigma}[f] \stackrel{def}{=} \sup_{\delta > 0}
\left[ \frac{\omega_{\sigma}[f](\delta)}{\delta^{\alpha}} \right].
$$
 \ Evidently, the  functional $ f \to H_{\alpha,\sigma}[f] $ is (complete) semi - norm relative the function
$ \ f, \ f \in HDT(\sigma, \alpha), $  as in the case of classical H\"older's norm, in which $ \sigma = 1. $\par

\vspace{3mm}

{\bf Example 4.1.} It is easily to compute by means of theorem 3.1,  that if $  f \in HDT(\sigma, \alpha), $ then

$$
\Delta_n[f] \le H_{\alpha, \sigma}[f] \cdot n^{-\alpha/2} \cdot \int_0^{\infty} z^{\alpha} \ | \ d Q_{\zeta}(z) \ |=
$$

$$
\alpha \cdot  H_{\alpha, \sigma}[f] \cdot n^{-\alpha/2} \cdot \int_0^{\infty} z^{\alpha - 1} \   Q_{\zeta}(z)  \ dz. \eqno(4.2)
$$

 \ If in addition the r.v. $ \zeta $ satisfies the condition of the example 3.1, then

$$
\Delta_n[f] \le K^{-\alpha/q} (p) \cdot \alpha \cdot \ H_{\alpha, \sigma}[f] \cdot n^{-\alpha/2} \cdot
 \Gamma \left( \alpha /q \right). \eqno(4.3)
$$

\vspace{3mm}

 \ The case $  p = q = 2 $  and $  K(p) = 1/2 $ correspondent to the classical Bernstein's case,  see \ \cite{Ostrovsky201}. \par

\vspace{4mm}

{\bf Remark 4.1.} The case when $ \omega_{\sigma}[f](\delta) $ is (continuous) non-negative {\it regular varying }
at the origin function:

$$
\omega_{\sigma}[f](\delta) \le H_L \cdot \delta^{\alpha} \ L(\delta), \ \delta \ge 0, \ \alpha = \const \in (0,1]
$$
where $ L = L(\delta)  $ is non-negative continuous in the set $ (0, b - a) $ slowly varying at the origin function,
that is

$$
\forall z > 0 \ \Rightarrow \lim_{\delta \to 0+} \frac{L(\delta z)}{L(\delta)} = 1,
$$
may be considered analogously. Indeed:

$$
\Delta_n[f] \le H_L \cdot \int_0^{\infty} \ L \left(\frac{z}{\sqrt{n}} \right) \ \left[ \frac{z}{\sqrt{n}} \right]^{\alpha} \ |d Q_{\zeta}(z)| \sim
$$

$$
H_L \cdot n^{-\alpha/2} \ L \left( \frac{1}{\sqrt{n}} \right) \int_0^{\infty} \ z^{\alpha} \ |d Q_{\zeta}(z)|, \ n \to \infty.
$$

\vspace{4mm}

 {\bf  Example 4.2;  "Poisson" case.} Suppose now that the r.v. $  \xi $ has a Poisson distribution with a parameter $ x. $
Here $   I = [0,\infty),  $ and let $ x = \lambda  \ge 1. $ Then

$$
{\bf P}(\xi = k) = e^{-x} \frac{x^k}{k!}, \ k = 0,1, \ldots.
$$

 \ Recall that "Poisson" approximation  of an uniform continuous function $  f  $ has a following form

$$
S_{n,P} = S_n[f](x) = e^{-nx} \sum_{k=0}^{\infty} \frac{(n x)^k}{k!} \ f \left( \frac{k}{n}  \right),
$$
and was investigated,  e.g. in the articles \cite{Chlodovsky1}, \cite{Mursalem1}, \cite{Szasz1}.\par

 \ Let us estimate the moment generating function for the centered and normed variable

$$
\eta = \frac{\eta_0}{\sqrt{x} }= \frac{\xi - x}{\sqrt{x}}.
$$

  \ We have for the values  $  z > 0 $ and $ \ \lambda = x \ge 1: $

$$
{\bf E} e^{\xi \lambda} = \sum_{k=0}^{\infty} e^{kz} e^{-\lambda} \frac{\lambda^k}{k!}  =
e^{-\lambda } \sum_{k=0}^{\infty} \frac{ (\lambda \ e^z )^k}{k!} = \exp \left( \lambda (e^z - 1) \right);
$$

$$
\ln {\bf E} e^{ z \eta} = -z \sqrt{\lambda} + \lambda \left(e^{ z /\sqrt{\lambda} } - 1 \right) =
$$

$$
\frac{z^2}{2} + \frac{z^3}{3! \ \lambda^{1/2} }  + \ldots + \frac{z^k}{k! \ \lambda^{ k/2 - 1 }} + \ldots \ . \eqno(4.4)
$$

 \ It is clear that the right-hand side  of the relation (4.4) attains the maximal value relative the variable $  \lambda $
at the point  $ \lambda = 1. $ Thus,

$$
\phi(z) = \phi_P(x) \stackrel{def}{=}  \sup_{\lambda \ge 1} {\bf E} \ e^{ z \eta} = e^z - 1 - z. \eqno(4.5)
$$

  \ The Young-Fenchel transform of these function has a form

$$
\ln \phi_P^*(u)  = u \ \ln(1 + u) - u + \ln (1 + u), \ u > 0,
$$
following $ \ln \phi_P^*(u)  \sim u \ \ln(1 + u), \ u \to \infty $ and analogously

$$
Q'_P(u) \sim \left\{ \ln(1 + u) + u/(u + 1) \right\} \cdot e^{-u \ \ln(1 + u)}, \ u \to \infty,
$$
the index "P" correspondent the name "Poisson." It remains to use theorem 3.1:

$$
\Delta_{n,P}[f] \le \int_0^{\infty}  \omega_{\sigma} [f] \left( \frac{z}{\sqrt{n}}  \right) \  Q'_{P}(z) \ dz. \eqno(4.6)
$$

 \  Here $ \sigma = \sigma(x) = \sqrt{x}, \ x \ge 1. $ \par

 \ Of course, if the function $  f(\cdot) $ satisfies H\"older's, more exactly, $ HDT $ condition (4.1), then the
integral in the right-hand side of inequality (4.6) converges and we conclude as before

$$
\Delta_{n,P}[f] \le  C_P \cdot H_{\alpha, \sigma}[f] \cdot n^{-\alpha/2}, \ n \ge 1.
\eqno(4.7)
$$

\vspace{3mm}

\section{ Low bounds. Exactness of our estimates.}

\vspace{3mm}

 \hspace{5mm} Let us consider the following example.  $  I := [0,1], $ the distribution $  \xi $ is such that $ \sigma (x) = \Var \xi  $
is continuous  and denote $ \overline{\sigma} = \max \sigma(x), \ x_0 = \argmax \sigma(x) \in (0,1). $ \par
 \ Note that the Bernstein's case  $  \sigma(x) = \sqrt{x(1 - x)}, \ x_0 = 1/2 $   is suitable for us. \par

 \ Let also $  g = g(x), \ x \in [0,1] $ be non-negative trial continuous function from the set $  HDT(\sigma,\alpha), \ 0 < \alpha \le 1 $
for which $  g(x_0) = 0 $ and $  |g(x_0 + \delta) - g(x_0)| =   $

$$
g(x_0 + \delta) \ge H_{\sigma, \alpha}[g] \cdot \delta^{\alpha} = H_- \cdot \delta^{\alpha},
 \ \delta \in [ 0,\min(x_0, 1 - x_0)]. \eqno(5.1)
$$

 \ We deduce

$$
\Delta_n[g] \ge {\bf E} g(x_0 + \zeta_n \overline{\sigma}/\sqrt{n} )  \ge  H_- \cdot n^{-\alpha/2} \cdot {\bf E} |\zeta_n|^{\alpha}. \eqno(5.2)
$$
 \ It follows from the (another) Bernstein's theorem that

$$
\lim_{n \to \infty}  {\bf E} |\zeta_n|^{\alpha} = {\bf E} |\tau|^{\alpha} = 2^{\alpha/2} \pi^{-1/2} \Gamma ( (\alpha + 1)/2  )
\stackrel{def}{=} G(\alpha),
$$
where the r.v. $  \tau $ has a standard normal distribution. \par

 \ To summarize, we have proved in fact the following {\it low} bound the generalized Bernstein's approximation. \par

\vspace{3mm}

{\bf Proposition 5.1.}

$$
\sup_{\const \ne g \in H(\alpha,\sigma)} \underline{\lim}_{n \to \infty}
\left[ \frac{\Delta_n[g]}{H_{\alpha,\sigma}[g] \ n^{-\alpha/2}}  \right]
\ge G(\alpha). \eqno(5.3)
$$

 \ See also \cite{Tikhonov1}.\par

 \vspace{3mm}

 \section{ Concluding remarks.}

 \vspace{3mm}

 \hspace{4mm} {\bf A.}  \ It is no hard perhaps to generalize by our opinion obtained results into the "more" multivariate case
$  d = 2,3,4,5, \ldots $ as well as into other methods of approximation, if only they had a probabilistic representation.\par

\vspace{3mm}

{\bf B.} \ One can also investigate and improve the rate of convergence of partial derivatives for the multivariate
 Bernstein's polynomials, in the spirit of the article \cite{Veretennikov1}. \par


\vspace{4mm}

\begin{thebibliography}{99}

\vspace{3mm}

\bibitem{Achiezer1}
{\sc  Achiezer (Akhiezer) N.I.} {\it Theory of approximation.}  Frederick Ungar Publishing Co., New York,  1956.

\bibitem{Bernstein1}
{\sc Bernstein S.N.} {\it Demonstration du theoreme de Weierstrass, fonde sur le probabilites. } Comm.
Sot. Math. Kharkow, {\bf 13}, (1912-1913), l-2.

\bibitem{Bojanic1}
{\sc Bojanic Ranko.} {\it On the approximation of continuous functions by Bernstein
polynomials.} Acad. Serbe Sci. Arts Glus, 232, (1959), 59-65. [In Serbian].

\bibitem{Bojanic2}
{\sc Bojanic Ranko and Fuhua Cheng.} {\it Rate of Convergence of Bernstein Polynomials
for Functions with Derivatives of Bounded Variation.} Journal of Mathematical
Analysis and Applications. {\bf 141,} \ 136-151, (1989).

\bibitem{Buldygin1}
{\sc Buldygin V.V., Kozatchenko Yu.V.} {\it About subgaussian  random variables. } Ukrainian Math. Journal, 1980,
{\bf 32,} $ N^o $ 6, 723-730.

\bibitem{Buldygin2}
{\sc Buldygin V.V., Moskvichova K.K.} {\it The sub-Gaussian norm of a binary random variable.}
Theor. Probability and Math. Statist.
Vip. 86,  No. 86, 2013, Pages 33-49.

\bibitem{Buldygin3}
{\sc Buldygin V.V., Kozatchenko Yu.V.} {\it Metric Characterization of Random
Variables and Random Processes.} 1998, Translations of Mathematics Monograph, AMS, v.188.

\bibitem{Chen1}
{\sc S. X. Chen and J. S. Liu.} {\it Statistical applications of the Poisson-binomial and conditional Bernoulli
distributions.} Statist. Sinica, 7(4): 875-892, 1997.

\bibitem{Chlodovsky1}
{\sc Chlodovsky, I. N.} {\it On some properties of S. N. Bernstein polynomials.} Proc.
1st All-Union congress of mathematics (Kharkov, 1930), Moscow; Leningrad:
ONTI NKTP USSR, 1936, {\it 22}, (in Russian).

\bibitem{Chlodovsky2}
{\sc I. Chlodovsky.} {\it Sur le d'eveloppement des fonctions d'efinies dans un intervalle infini en s'eries de
polynomes de S.N. Bernstein.}  Compos. Math. 4 (1937) 380-393.

\bibitem{DeVore1}
{\sc DeVore, Lorentz G.G.} {\it Constructive approximation.} Springer Verlag, 1997.

\bibitem{Ditzian1}
{\sc Z. Ditzian, V. Totik.} {\it Moduli of Smoothness.} Springer, New York, 1987.

\bibitem{Floater1}
{\sc Floater M. S.} {\it On the convergence of derivatives of Bernstein approximation.}
Journal of Approximation Theory, 2005, 134(1), 130–135.

\bibitem{Gonska1}
{\sc  Gonska H., Zhou D-x. } {\it On an extremal problem concerning Bernstein
operators,} Serdica Math. J. 21, 1995, 137-150.

\bibitem{Glazyrina1}
{\sc Polina Glazyrina and Sergey Tikhonov.} {\it Jacobi weights, fractional integration,
and sharp Ulyanov inequalities.} \\
arXiv:1601.00814v1 [math.FA] 5 Jan 2016

\bibitem{Kahane1}
{\sc Kahane J.P.} {\it  Properties locales des fonctions a series de Fourier aleatoires.} Studia Math. (1960), \  {\bf 19,}
$  N^o $ 1, 1-25.

\bibitem{Kozatchenko1}
{\sc  Kozatchenko Yu.V., Ostrovsky E.I.} {\it Banach spaces of random
variables of subgaussian type.} Theory Probab. And Math. Stat.,
Kiev, (1985), p. 42 $ \- \ $ 56 (in Russian).

\bibitem{Li2}
{\sc  Li Z.} {\it Bernstein polynomials and modulus of continuity.} J. Approx. Theory, {\bf 102,}
(2000),  171-174; MR1736050 (2000j:41029).

\bibitem{Lorentz1}
{\sc Lorentz G.G.} {\it Bernstein Polynomials.} Univ. of Toronto Press, Toronto, 1953.

\bibitem{Math'e}
{\sc  Math'e P. }  {\it Approximation of H\"older continuous functions by Bernstein polynomials.}
Amer. Math. Monthly, {\bf 106}, (1999), 568-574.

\bibitem{Mursalem1}
{\sc M. Mursaleen and Khursheed J. Ansari.}  {\it Approximation by generalized Sz'asz operators
involving Sheffer polynomials.} \\
arXiv:1601.00675v1 [math.CA] 25 Dec 2015

\bibitem{Ostrovsky1}
{\sc Ostrovsky E.I. } {\it Exponential Estimations for Random Fields.} Monograph,
Moscow-Obninsk, OINPE, (1999), (in Russian).

\bibitem{Ostrovsky102}
{\sc Ostrovsky E., Sirota L.} {\it  Subgaussian and strictly subgaussian random variable. }
arXiv:1406.3933v1 [math.PR] 16 Jun 2014


\bibitem{Ostrovsky103}
{\sc Ostrovsky E.I. } {\it  Generalization of the norm Buldygin - Kosatchenko and Central Limit Theorem in Banach space.}
 Probab. Theory Appl., (1982), issue 3, p. 616-618. (in Russian).

\bibitem{Ostrovsky201}
{\sc Ostrovsky E., Sirota L.} {\it Non-uniforem non-asymptotical sharp estimate for the rate of convergence for
Bernstein's polynomial approximation, with bilateral constant evaluation. } \\
arXiv:1508.06933v1 [math.FA] 27 Aug 2015

\bibitem{Palatanea1}
{\sc P\"alat\"anea Radu.} {\it On some constants in approximation by Bernstein operators.}
General Mathematics,  Vol. 16, No. 4, (2008), 137-148.

\bibitem{Palatanea2}
{\sc P\"alat\"anea Radu.} {\it  Optimal constant in approximation by Bernstein operators. }
J. Comput. Analysis Appl., 5 (2), 2003, 195-235.

\bibitem{Stancu1}
{\sc D.D. Stancu. } {\it Approximation of functions by means of a new generalized Bernstein operator.}
 Calcolo, {\bf 20,} (2,) (1983),  211-229.

\bibitem{Szasz1}
{\sc O. Sz'asz.} {\it Generalization of S. Bernstein's polynomials to the infinite interval.} J. Research Nat.
Bur. Standards, {\bf 45,} (1950),  239-245.

\bibitem{Telyakovskii1}
{\sc Telyakovskii S. A.} {\it On the rate of approximation of functions by the Bernstein
polynomials.} Proc. Inst. Math. Mech., 2009, 264, suppl. 1, pp. 177-184.

\bibitem{Tikhonov1}
{\sc Tikhonov, I. V., Sherstyukov, V. B.}  {\it The module function approximation by
Bernstein polynomials.} Cheliabinsk State Univ. Vestnik, Ser. Matem., Mech.,
Inform., 2012, {\bf 15,} (26), 6-40. (in Russian).

\bibitem{Timan1}
{\sc Timan A.F.} {\it Theory of Approximation of Functions of a Real Variable.} International Series of
Monographs on Pure and Applied Mathematics, Paperback–January 1, 1963. Pergamon Press.

\bibitem{Veretennikov1}
{\sc  Veretennikov A.Yu., Veretennikova E.V.}
{\it On convergence of partial derivatives of multivariate Bernstein polynomials.}\\
arXiv:1507.05235v1 [math.PR] 18 Jul 2015

\vspace{4mm}

\end{thebibliography}
\end{document}